\documentclass[reqno]{amsart}
\usepackage[utf8]{inputenc}
\usepackage[T1]{fontenc}
\usepackage{amsthm}
\usepackage{amsmath}
\usepackage{amssymb}
\usepackage{xcolor}
\usepackage{mathtools} 
\usepackage{color}
\usepackage{hyperref} 
\usepackage{xfrac}
\usepackage{todonotes}
\usepackage[shortlabels]{enumitem}
\hypersetup{
	colorlinks,
	linkcolor={red!50!black},
	citecolor={blue!30!black},
	urlcolor={blue!30!black}
}
\usepackage[nameinlink, capitalise, noabbrev]{cleveref}

\newtheorem{theorem}{Theorem}

\newtheorem{kor}[theorem]{Corollary}

\theoremstyle{definition}

\crefname{kor}{Corollary}{Corollaries}
\crefname{lem}{Lemma}{Lemmata}

\numberwithin{equation}{section}
\numberwithin{theorem}{section}
\numberwithin{definition}{section}
\DeclareMathOperator\sgn{sgn}

\title{A new Berry-Esseen-type estimate in the \\ free central limit theorem}
\author{Leonie Neufeld}
	\thanks{Fakult\"at f\"ur Mathematik,
	Universit\"at Bielefeld, 33501 Bielefeld, 
	Germany;  lneufeld@math.uni-bielefeld.de}
\thanks{Funded by the Deutsche Forschungsgemeinschaft (DFG, German Research Foundation) -- IRTG 2235 -- 282638148.}
\begin{document} 
	
\begin{abstract}
Using the subordination approach, we provide a new Berry-Esseen-type estimate in the free central limit theorem in terms of the fourth Lyapunov fraction. In the special case of identical distributions, our result implies a rate of order $n^{-1/2 + \varepsilon}$ for any $\varepsilon>0$, thus almost leading to the optimal rate of order $n^{-1/2}$. 
\end{abstract}	
	\maketitle
\section{Introduction} \label{Intro}
Let $(X_i)_{i \in \mathbb{N}}$ be a sequence of free self-adjoint random variables with analytic distributions $(\mu_{X_i})_{i \in \mathbb{N}}$. Among other moment conditions, we assume that each $\mu_{X_i}$ has mean zero and variance $\sigma_i^2 \in (0, \infty)$. The aim of this work is to analyze at what rate the analytic distribution $\mu_{S_n}$ of the normalized sum $S_n$ defined by
\begin{align*}
S_n := \frac{1}{B_n} \sum_{i=1}^{n} X_i, \qquad	B_n := \left( \sum_{i=1}^{n} \sigma_i^2 \right)^{\frac{1}{2}},
\end{align*}
converges weakly to the Wigner semicircle distribution $\omega$. 

\subsection{Known rates of convergence in the free CLT}
Before we state our result, let us briefly recall what is already known about the asymptotic behavior of the analytic distribution $\mu_{S_n}.$ 

The \textit{free central limit theorem (free CLT)}, initially proved by Voiculescu \cite{Voiculescu1985} and further extended in \cite{Speicher1990,Maassen1992,Pata1996,Kargin2007b,Maejima2023}, provides conditions (such as the \textit{Lindeberg condition}) under which $\mu_{S_n}$ converges weakly to the (standard) Wigner semicircle distribution $\omega$. The associated rate of convergence has been studied in several papers.  In the following, we restrict to rates measured in terms of the \textit{Kolmogorov distance} $\Delta$ given by 
\begin{align*}
		\Delta(\mu, \nu) := \sup_{x \in \mathbb{R}} \big \vert \mu((-\infty,x]) - \nu((-\infty,x]) \big\vert
\end{align*}	
for any probability measures $\mu$ and $\nu$ on $\mathbb{R}.$ We refer to  \cite{Chistyakov2013,Chistyakov2017,Fathi2017, Cebron2020,Austern2020,Diez2023,Diaz2024} for rates of convergence in the free CLT with respect to other distances.

%and to \cite{Petrov1975} for rates in the setting of the classical CLT.

Under the usual finite third absolute moment assumption, and by using the concept of subordination,  Chistyakov and Götze \cite{Chistyakov2008,Chistyakov2013} established the optimal rate of order $n^{-\frac{1}{2}}$ (in $\Delta$) for the case of identical distributions, thus obtaining the same rate as in the classical CLT; see \cite[Chapter 5]{Petrov1975}. In the setting of non-identical distributions, they provided a rate of convergence that, up to a square root, corresponds to the classical rate given by the \textit{third Lyapunov fraction} $L_{3n}$. More precisely, denoting the finite third absolute moment of $\mu_{X_i}$ by $\beta_3(\mu_{X_i})$, Chistyakov and Götze proved
\begin{align} \label{Rate Goetze}
\Delta(\mu_{S_n}, \omega) \leq c L_{3n}^{\frac{1}{2}}, \qquad 	L_{3n} := \frac{\sum_{i=1}^n \beta_3(\mu_{X_i})}{B_n^3},
\end{align}
for some constant $c>0$.

Restricting to bounded random variables, a similar result was derived in \cite{Neufeld2024a}: If the analytic distribution $\mu_{X_i}$ has support in $[-L_i, L_i]$ for some $L_i >0$ and any $i \in \mathbb{N}$, then we obtain
\begin{align} \label{rate meins}
	\Delta(\mu_{S_n}, \omega)  \leq c L_{S, 3n}, \qquad	L_{S, 3n} := \frac{\sum_{i=1}^n L_i^3}{B_n^3},
\end{align}
for some constant $c>0.$ In particular, it is possible to remove the square root in \eqref{Rate Goetze} at the cost of replacing the Lyapunov fraction $L_{3n}$ by the greater Lyapunov-type fraction $L_{S,3n}$. Note that in the setting of identically distributed bounded random variables, the result in \eqref{rate meins} implies the optimal rate of order $\smash{n^{-\frac{1}{2}}}$, while \eqref{Rate Goetze} yields a rate of order $\smash{n^{-\frac{1}{4}}}$.

For more rates in the free CLT with respect to the Kolmogorov distance, proved by alternative methods or leading to faster convergence under vanishing \textit{free cumulants}, we refer to \cite{Kargin2007a,Salazar2023}.

%We refer to \cite{Kargin2007a,Austern2020} for more rates in terms of the Kolmogorov distance, and to \cite{Fathi2017, Chistyakov2017, Austern2020, Diaz2024} for rates of convergence in the free CLT with respect to other distances. 

\subsection{A new Berry-Esseen-type estimate in the free CLT} 
Let us now formulate the main theorem of this work, providing a new rate of convergence in the free CLT for not necessarily identically distributed and possibly unbounded random variables in terms of the associated \textit{fourth Lyapunov fraction} $L_{4n}$.

\begin{theorem} \label{main theorem}
 Let $(X_i)_{i \in \mathbb{N}}$ be a sequence of free self-adjoint random variables with analytic distributions $(\mu_{X_i})_{i \in \mathbb{N}}$. Assume that each $\mu_{X_i}$ has mean zero, variance $\sigma_i^2>0$, and finite fourth moment $m_4(\mu_{X_i})$. Define 
\begin{align*}
S_n := \frac{1}{B_n} \sum_{i=1}^n X_i, \qquad B_n :=\left( \sum_{i=1}^n \sigma_i^2 \right)^{\frac{1}{2}}, 
\end{align*}
and let $\mu_{S_n}$ denote the analytic distribution of $S_n$.
Then, for any $\varepsilon \in (0, \frac{1}{2})$, we have 
\begin{align*}
\Delta(\mu_{S_n}, \omega) \leq C_\varepsilon L_{4n}^{\frac{1}{2} - \varepsilon}, \qquad L_{4n}:= \frac{\sum_{i=1}^{n} m_4(\mu_{X_i})}{B_n^4},
\end{align*}
for some constant $C_\varepsilon>0$.  
\end{theorem}

The proof of the above theorem has a recursive structure and is based on modifications of the subordination method from \cite{Chistyakov2008,Chistyakov2013,Neufeld2024a}. As can be derived from the proof (more precisely from \cref{corollary}), the constant $C_\varepsilon$ depends on the parameter $\varepsilon$ via the number of recursions.

%The dependency of the constant $C_\varepsilon$ on $\varepsilon$ will be evident from \cref{corollary}. Instead, 
Finally, let us compare the rate in \cref{main theorem} with those in \eqref{Rate Goetze} and \eqref{rate meins}: 
There are cases where \cref{main theorem} implies a better rate of convergence, and examples in which \eqref{Rate Goetze} or \eqref{rate meins} provide better rates. However, similar to the result in \eqref{rate meins}, the above theorem yields almost the optimal rate of order $n^{-\frac{1}{2}}$ in the case of identical distributions. In more detail, if the random variables $(X_i)_{i \in \mathbb{N}}$ from \cref{main theorem} are identically distributed, then the associated fourth Lyapunov fraction $L_{4n}$ is of order $n^{-1}$, which leads to
\begin{align*}
	\Delta(\mu_{S_n}, \omega) \leq \left(\frac{m_4(\mu_{X_1})}{\sigma_1^4}\right)^{\frac{1}{2}-\varepsilon} \frac{C_\varepsilon}{n^{\frac{1}{2}-\varepsilon}}
\end{align*}
for any $\varepsilon \in (0, \frac{1}{2}).$ 

%\subsection*{Organization} In Section \ref{Section: Preliminaries} we briefly recall a few basics of free probability theory. Section \ref{Section: Proof} is devoted to the proof of Theorem \ref{main theorem}.

\subsection*{Acknowledgments} I would like to thank Friedrich Götze for valuable discussions.

\pagebreak
\section{Preliminaries} \label{Section: Preliminaries}
In this section, we fix some notation and recall the subordination approach to free additive convolutions as well as Bai's smoothing inequality. 

\subsection{Notation} \label{Pre: Notation}
Throughout this work, we let $[n] := \{1, \dots, n\}$ for any $n \in \mathbb{N}$.

Moreover, when writing $\sqrt{z}$ for some $z \in \mathbb{C} \setminus [0, \infty)$, we refer to the complex square root with branch cut placed on the non-negative real axis. Note that we have
\begin{align} \label{square root real and im formula}
	\Re\sqrt{z} = \sgn( \Im z)\sqrt{\frac{1}{2} \left(\sqrt{(\Re z )^2+(\Im z)^2} + \Re z\right)}, \qquad \Im\sqrt{z} =  \sqrt{\frac{1}{2}\left( \sqrt{(\Re z)^2+(\Im z)^2} - \Re z\right)} \geq 0
\end{align}
for $z  \in \mathbb{C}\setminus[0, \infty)$. Above, $\sgn$ denotes the sign function under the convention $\sgn(0)=1$. 

Finally, all probability measures in this paper are assumed to be Borel probability measures on $\mathbb{R}$. For such a probability measure $\mu$, we let
\begin{align*}
	m_k(\mu) := \int_{\mathbb{R}} x^k \mu(dx), \qquad \beta_k(\mu) := \int_{\mathbb{R}} \vert x \vert ^k \mu(dx), \qquad k \in \mathbb{N},
\end{align*}
denote its $k$-th (absolute) moments.

\subsection{Subordination approach to free additive convolutions} \label{subordination approach}
As explained in \cref{Intro}, the aim of this work is to study the analytic distribution of the sum of certain free self-adjoint random variables. Such distributions are also known as \textit{free additive convolutions}. More precisely, given random variables $X_1, \dots, X_n$ as above with analytic distributions $\mu_{X_1}, \dots, \mu_{X_n}$, the analytic distribution of $X_1 + \dots +  X_n$ is called the free additive convolution of $\mu_{X_1}, \dots, \mu_{X_n}$, and is denoted by $\mu_{X_1} \boxplus \dots \boxplus \mu_{X_n}$. Using the \textit{subordination approach} (see \cref{subordination equations} below), one can define the free additive convolution as an operation on the set of probability measures on the real line, without any reference to random variables and their analytic distributions. We refer to  \cite{Voiculescu1985, Voiculescu1986, Maassen1992, Bercovici1993} for extensive accounts on free additive convolutions, and to \cite{Belinschi2007,Chistyakov2011} for details on the subordination approach.

Let $\mu$ be a probability measure on $\mathbb{R}$. The Cauchy transform $G_\mu$ and the $F$-transform $F_\mu$ of $\mu$ are defined by
\begin{align*}
	G_\mu(z) := \int_{\mathbb{R}} \frac{1}{z-x} \mu(dx), \qquad	F_{\mu}(z) := \frac{1}{G_{\mu}(z)}, \qquad z \in \mathbb{C}^+,
\end{align*}
where $\mathbb{C}^+$ denotes the complex upper half-plane. For later reference, we note that the corresponding transforms of the Wigner semicircle distribution $\omega$ are given by 
\begin{align*}
	G_\omega(z) = \frac{1}{2}\left( z - \sqrt{z^2-4} \right), \qquad	F_\omega(z) = \frac{1}{2}\left( z + \sqrt{z^2-4} \right), \qquad z \in \mathbb{C}^+,
\end{align*}
and that we have
\begin{align} \label{bound CT omega 1}
	\vert G_\omega(z) \vert \leq 1, \qquad z \in \mathbb{C}^+; 
\end{align}
see \cite[Lemma 8]{Kargin2007a}.

The subordination approach to free additive convolutions can be formulated on the basis of the $F$-transform, as demonstrated in the following theorem. We refer to \cite[Corollary 2.2]{Chistyakov2011} for a proof.
% We remark that this approach allows to define free additive convolution of probability measures on $\mathbb{R}$ solely by the use of the corresponding reciprocal Cauchy transforms, without any reference to non-commutative random variables. 

\begin{theorem}\label{subordination equations}
	Let $\mu_1, \dots, \mu_n$ be probability measures on $\mathbb{R}$. There exist unique holomorphic functions $Z_1, \dots, Z_n: \mathbb{C}^+ \rightarrow \mathbb{C}^+$ such that for any $z \in \mathbb{C}^+$ the equations
	\begin{align*}
		\left(\sum_{i=1}^n Z_i(z)\right)- z = (n-1)F_{\mu_1}(Z_1(z)), \qquad F_{\mu_1}(Z_1(z)) = \cdots = F_{\mu_n}(Z_n(z)) 
	\end{align*}
	hold. The so-called subordination functions $Z_1, \dots, Z_n$ satisfy $\Im Z_i(z) \geq \Im z$ for all $z \in \mathbb{C}^+$  and $ i \in [n].$ Moreover, there exists a probability measure $\mu$ on $\mathbb{R}$ with $F_\mu(z) = F_{\mu_1}(Z_1(z))$ for all $z \in \mathbb{C}^+$. We set $\mu_1 \boxplus \cdots \boxplus \mu_n := \mu$. 
\end{theorem}

\subsection{Bai's smoothing inequality}
Recall that the Kolmogorov distance $\Delta(\mu, \nu)$ between two probability measures $\mu$ and $\nu$ on $\mathbb{R}$ is defined by 
\begin{align*}
	\Delta(\mu, \nu) := \sup_{x \in \mathbb{R}} \big \vert \mu((-\infty,x]) - \nu((-\infty,x]) \big\vert.
\end{align*}			

The following theorem, which is a variation of a well-known smoothing inequality due to Bai \cite{Bai1993}, provides an upper bound on the Kolmogorov distance between a probability measure and the Wigner semicircle distribution in terms of their Cauchy transforms. For a proof, we refer to \cite[Section 2]{Goetze2003}.

\begin{theorem} \label{Bai}
	Let $\mu$ be a probability measure on $\mathbb{R}$ with 
	\begin{align*} 
		\int_{-\infty}^\infty \big \vert \mu((-\infty, x]) - \omega((-\infty, x]) \big \vert dx < \infty.
	\end{align*}
Choose $v \in (0,1)$, $\varepsilon \in (0,2)$, and $a, \gamma >0$ in such a way that
\begin{align*}
	\gamma = \frac{1}{\pi} \int_{\vert x \vert < a} \frac{1}{1+x^2} dx > \frac{1}{2} \qquad \text{and} \qquad \varepsilon>2va
\end{align*}
are satisfied. Define $I_\varepsilon := [-2+\frac{\varepsilon}{2}, 2- \frac{\varepsilon}{2}]$. Then, we have 
	\begin{align*}
	\Delta(\mu, \omega) & \, \leq C_{\gamma} \left( \frac{4a^2 v}{\pi} + \gamma \varepsilon^{\frac{3}{2}}  + \int_{-\infty}^{\infty} \left\vert G_\mu(u+i) - G_\omega(u+i) \right\vert du  +  \sup_{x \in I_\varepsilon} \int_{v}^1 \vert G_\mu(x+iy) - G_\omega(x+iy) \vert dy  \right),
\end{align*}
where $C_\gamma>0$ is given by $C_\gamma := ((2\gamma -1)\pi)^{-1}$.
\end{theorem}

\section{Proof of  \texorpdfstring{\cref{main theorem}}{Theorem 1.1}}  \label{Section: Proof}
In this section, we prove \cref{main theorem}. The procedure is as follows: In Section \ref{Section: Preliminary}, we formulate  an auxiliary result (see \cref{main theorem 2}) and show how it implies \cref{main theorem}. Then, in \cref{Section: Preliminary Proof}, we verify \cref{main theorem 2}.
\subsection{An auxiliary result} \label{Section: Preliminary}
 \cref{main theorem} will follow from the following auxiliary result:
\begin{theorem} \label{main theorem 2}
	Let $S_n, \mu_{S_n}$, and $L_{4n}$ be as in \cref{main theorem}.
	Assume that there exist $k \in [0, \frac{1}{2})$ and a constant $D(k)>0$  such that $\smash{\Delta(\mu_{S_n}, \omega) \leq D(k)L_{4n}^k}$
	holds. Then, we have
	\begin{align*}
		\Delta(\mu_{S_n}, \omega) \leq C(k) L_{4n}^{\frac{1}{4} + \frac{k}{2}}
	\end{align*}
	for some constant $C(k)>0$, which can be chosen to depend on $k$ only via $D(k)$.
\end{theorem}
Note that in the relevant case $L_{4n} <1$, the above theorem provides an improved rate of convergence: Due to $k < \frac{1}{2}$, we have $\frac{1}{4} + \frac{k}{2} >k$, which leads to 
\begin{align*}
L_{4n}^{\frac{1}{4} + \frac{k}{2}} < L_{4n}^k.
\end{align*}

It is evident that \cref{main theorem 2} allows for repeated applications. In particular, by induction, we obtain: 
\begin{kor} \label{corollary}
	Let $S_n, \mu_{S_n}$, and $L_{4n}$ be as in \cref{main theorem}.
Assume that there exist $k \in [0, \frac{1}{2})$ and a constant $D(k)>0$  such that $	\smash{\Delta(\mu_{S_n}, \omega) \leq D(k) L_{4n}^k}$
holds. Then, for any $i \in \mathbb{N}$ with $i \geq 1$,  we have
\begin{align*}
	\Delta(\mu_{S_n}, \omega) \leq C(i,k) L_{4n}^{\frac{2^i-1}{2^{i+1}} +  \frac{k}{2^i}}
\end{align*}
for some constant $C(i,k) >0$, which can be chosen to depend on $k$ only via $D(k)$.
 \end{kor}

\pagebreak
Finally, we show how the above corollary implies \cref{main theorem}: 

\begin{proof}[Proof of \cref{main theorem}] 
	Let $S_n, \mu_{S_n}$, and $L_{4n}$ be as in \cref{main theorem}. We consider two cases: 
	First, assume that $L_{4n}\geq 1$ holds. Since the Kolmogorov distance is bounded by $1$, we immediately obtain 
	\begin{align*}
	 \Delta(\mu_{S_n}, \omega) \leq 1 \leq L_{4n}^{\frac{1}{2}-\varepsilon}
	\end{align*}
	for any $\varepsilon \in (0, \frac{1}{2}).$
	
	Second, let $L_{4n}<1$. Combining \eqref{Rate Goetze} with the  inequality 
	\begin{align} \label{inequality lyapunov}
		L_{3n} \leq L_{4n}^{\frac{1}{2}},
	\end{align}
	see \cite[Chapter VI, §2, Lemma 2]{Petrov1975} for a proof, we observe that it suffices to restrict to $\varepsilon \in (0, \frac{1}{4})$. For such $\varepsilon$, define 
		\begin{align*}
		i(\varepsilon) := \left \lceil \frac{\log \left( \frac{1}{4\varepsilon} \right)}{\log (2)} \right \rceil \in \mathbb{N},
	\end{align*}
	where $\lceil \cdot \rceil$ denotes the ceiling function.
Then, we have 
\begin{align*}
	\frac{2^{i(\varepsilon)}-1}{2^{i(\varepsilon)+1}} + \frac{1}{4} \frac{1}{2^{i(\varepsilon)}}   \geq \frac{1}{2} - \varepsilon. 
\end{align*}
By \eqref{Rate Goetze} and \eqref{inequality lyapunov}, we deduce that the premise of \cref{corollary} is satisfied for $k=\frac{1}{4}$ and some constant $D(\frac{1}{4})>0$. 
Applying the last-named corollary to $i(\varepsilon)$ and $k$ as chosen above, we find a constant $C(i(\varepsilon),\frac{1}{4})>0$ such that
	\begin{align*}
	\Delta(\mu_{S_n}, \omega) \leq  C(i(\varepsilon),\tfrac{1}{4}) L_{4n}^{	\frac{2^{i(\varepsilon)}-1}{2^{i(\varepsilon)+1}} + \frac{1}{4} \frac{1}{2^{i(\varepsilon)}}} \leq C(i(\varepsilon), \tfrac{1}{4}) L_{4n}^{\frac{1}{2}-\varepsilon}
	\end{align*}
holds. Finally, by setting $C_\varepsilon := \max\big \{1, C(i(\varepsilon), \frac{1}{4}) \big \}$, the claim follows.
\end{proof}

\subsection{Proof of \texorpdfstring{\cref{main theorem 2}}{Theorem 3.1}} \label{Section: Preliminary Proof}
The proof of \cref{main theorem 2} consists of six steps: In the first step, we derive a quadratic functional equation for one of the subordination functions of the free additive convolution $\mu_{S_n}.$ Solving this equation, we obtain a precise formula for the subordination function in the second step. Then, in the third step, we apply Bai's inequality from \cref{Bai}, leaving us with two integrals to bound. In the fourth and fifth step, we establish suitable estimates for these integrals. We complete the proof of \cref{main theorem 2} in the sixth step. 

Before we begin with the first step, let us introduce some notation. Fix $n \in \mathbb{N}$. For $(X_i)_{i \in \mathbb{N}}$ and $B_n$ as in \cref{main theorem}, we let $\mu_i$ denote the analytic distribution of the normalized random variable $B_n^{-1}X_i$. Note that $\mu_{S_n} = \mu_1 \boxplus \dots \boxplus \mu_n$ holds. The subordination functions of $\mu_{S_n}$ are given by $Z_1, \dots, Z_n$. Let $G_1, \dots, G_n, G_{S_n}$ be the Cauchy transforms of $\mu_1, \dots, \mu_n, \mu_{S_n}$, while the corresponding $F$-transforms are xdenoted by $F_1, \dots, F_n, F_{S_n}$. Define $\rho_i^2 := m_2(\mu_i) = B_n^{-2}\sigma_i^2$ and observe that we have $\sum_{i=1}^n \rho_i^2=1$. Without loss of generality, assume that $\rho_1^2 = \min_{i \in [n]} \rho_i^2$ and $L_{4n}<1$ hold.

\subsubsection*{Step 1: Quadratic functional equation for $Z_1$} 
In this step, we derive a quadratic functional equation for the subordination function $Z_1.$
By \cref{subordination equations}, it follows
\begin{align*}
	Z_1(z) - z = \sum_{i=2}^{n} F_i(Z_i(z)) - Z_i(z), \qquad z \in \mathbb{C}^+,
\end{align*}
from which we deduce
\begin{align*}
Z_1(z) - z + \frac{1}{Z_1(z)} = \frac{M_1(z) + M_2(z) + M_3}{Z_1(z)} = \frac{q(z)}{Z_1(z)}
\end{align*}
with 
\begin{align*}
M_1(z) &:= Z_1(z) \left( \sum_{i=2}^{n} F_i(Z_i(z)) - Z_i(z)+ \frac{\rho_i^2}{Z_i(z)} \right), \qquad 
M_2(z) := Z_1(z) \left( \sum_{i=2}^{n} \frac{\rho_i^2}{Z_1(z)} - \frac{\rho_i^2}{Z_i(z)} \right), \\
 &\qquad \qquad \qquad \qquad \,\,\,\,\,\,\,  M_3 := \rho_1^2,\qquad q(z) := M_1(z) + M_2(z) + M_3
\end{align*}
for all $z \in \mathbb{C}^+.$ In particular, for
\begin{align*}
Q(z, \omega) := \omega^2 - z\omega + 1 - q(z),
\end{align*}
we obtain $Q(z, Z_1(z)) = 0$ for all $z$ as above. 

In what follows, our goal is to bound the terms $M_1(z)$, $M_2(z)$, and $M_3$ appropriately. For this, we need to do some preparatory work. \cref{subordination equations} and the definition of the Cauchy transform yield
\begin{align*}
	r_{n,i}(z) :=   G_{S_n}(z)Z_i(z) - 1  = G_i(Z_i(z))Z_i(z) - 1 = \int_{\mathbb{R}} \frac{u}{Z_i(z) - u} \mu_i(du), \qquad z \in \mathbb{C}^+, i \in [n].
\end{align*}
For later reference, we observe that 
\begin{align*}
F_i(Z_i(z)) - Z_i(z) = - \frac{r_{n,i}(z)}{ 1+ r_{n,i}(z)}Z_i(z)
\end{align*}
and 
\begin{align} \label{subordination times cauchy expansion}
	G_i(Z_i(z))Z_i(z)& = 1 + \frac{1}{Z_i(z)} \int_{\mathbb{R}} \frac{u^2}{Z_i(z) - u} \mu_i(du)  \nonumber\\
	& =  1 + \frac{\rho_i^2}{Z_i^2(z)} + \frac{1}{Z_i^2(z)} \int_{\mathbb{R}} \frac{u^3}{Z_i(z) - u} \mu_i(du)   \\
	&  = 1 +  \frac{\rho_i^2}{Z_i^2(z)} + \frac{m_3(\mu_i)}{Z_i^3(z)} +  \frac{1}{Z_i^3(z)}  \int_{\mathbb{R}} \frac{u^4}{Z_i(z) - u}\mu_i(du)  \nonumber
\end{align}
hold for all $ z \in \mathbb{C}^+$ and $i \in [n]$. Let us continue by bounding $r_{n,i}(z)$ for certain $z$. Using Cauchy's inequality, the identity 
\begin{align*}
\frac{1}{\vert Z_i(z) - u \vert^2} = - \Im \left( \frac{1}{Z_i(z)-u} \right) \frac{1}{\Im Z_i(z)}, \qquad z \in \mathbb{C}^+, i \in [n], u \in \mathbb{R},
\end{align*}
and \cref{subordination equations},
we get 
\begin{align*}
	\vert r_{n,i}(z) \vert  \leq \left( \int_{\mathbb{R}} u^2 \mu_i(du) \right)^{\frac{1}{2}} \left( \frac{1}{\Im Z_i(z)} \int_{\mathbb{R}} - \Im\left( \frac{1}{Z_i(z)-u}\right) \mu_i(du) \right)^{\frac{1}{2}} \leq \frac{ \vert \rho_i \vert }{\sqrt{\Im z}} \vert G_i(Z_i(z)) \vert^{\frac{1}{2}}  
\end{align*}
for all $z \in \mathbb{C}^+$ and $i \in [n].$
By \eqref{bound CT omega 1}, integration by parts, and the premise of \cref{main theorem 2}, it follows
\begin{align} \label{integration by parts G_Sn}
	\begin{split}
\vert G_i(Z_i(z)) \vert = \vert G_{S_n}(z) \vert &\leq \vert G_\omega(z) \vert + \vert G_{S_n}(z) - G_\omega(z) \vert \leq 1+  \vert G_{S_n}(z) - G_\omega(z) \vert \\ & \leq 1+ \frac{\pi \Delta(\mu_{S_n}, \omega)}{\Im z} \leq 1+ \frac{\pi D(k) L_{4n}^{k}}{\Im z}
	\end{split}
\end{align}
for $z$ and $i$ as before and $k \in [0, \tfrac{1}{2})$ from \cref{main theorem 2}.
Hence, together with $0< \vert \rho_i \vert \leq L_{4n}^{\frac{1}{4}}$ for all $i \in [n]$, we deduce
\begin{align*}
\vert r_{n,i}(z) \vert \leq \frac{ \vert \rho_i \vert}{\sqrt{\Im z}} \left( 1 + \frac{\pi D(k) L_{4n}^k}{\Im z} \right)^{\frac{1}{2}} \leq \frac{L_{4n}^{\frac{1}{4}}}{\sqrt{\Im z}} \left( 1 + \frac{\sqrt{\pi D(k)} L_{4n}^{\frac{k}{2}}}{\sqrt{\Im z}} \right) = \frac{L_{4n}^{\frac{1}{4}}}{\sqrt{\Im z}} + \frac{\sqrt{\pi D(k)} L_{4n}^{ \frac{1}{4}+ \frac{k}{2}}}{\Im z} 
\end{align*}
for $z \in \mathbb{C}^+$ and $i \in [n]$. Now, define 
\begin{align*}
	D_1 := \left \{ z \in \mathbb{C}^+: \Im z \geq C_1 L_{4n}^{ \frac{1}{4}+ \frac{k}{2}} \right \}
\end{align*}
for some constant $C_1 = C_1(k)>0$, which satisfies $\smash{C_1^{-\frac{1}{2}} + \sqrt{\pi D(k)}C_1^{-1} < \frac{1}{10}}$. Due to $k < \frac{1}{2}$ and $L_{4n} <1$, we obtain 
\begin{align} \label{bound r 0.1}
\vert r_{n,i}(z) \vert \leq \frac{L_{4n}^{\frac{1}{8}-\frac{k}{4}}}{\sqrt{C_1}} + \frac{\sqrt{\pi D(k)}}{C_1} <  \frac{1}{\sqrt{C_1}} + \frac{\sqrt{\pi D(k)}}{C_1} < \frac{1}{10}
\end{align}
for all $z \in D_1$ and $i \in [n].$  

We proceed by deriving a lower bound on the modulus of the subordination functions. 
Using the definition of $r_{n,i}(z)$ as well as \eqref{integration by parts G_Sn} and \eqref{bound r 0.1}, it follows
\begin{align} \label{lower bound subordination}
\begin{split}
 \vert Z_i(z) \vert = \left \vert \frac{r_{n,i}(z) + 1 }{G_{S_n}(z)}\right \vert &\geq \frac{1 - \vert r_{n,i}(z) \vert}{\vert G_{S_n}(z) \vert}  \\ & \geq  \frac{\frac{9}{10}}{1 + \frac{\pi D(k) L_{4n}^k}{\Im z}}  \geq \frac{\frac{9}{10}}{1 + \pi D(k) C_1^{-1} L_{4n}^{\frac{k}{2} - \frac{1}{4}}} \\ &> \frac{\frac{9}{10}}{(1 + \pi D(k)C_1^{-1}) L_{4n}^{\frac{k}{2} - \frac{1}{4}}} = C_2L_{4n}^{\frac{1}{4} - \frac{k}{2}}
\end{split} 
\end{align}
for all $z \in D_1$ and $i \in [n]$ with $C_2 = C_2(k) := 9 (10(1 + \pi D(k)C_1^{-1}))^{-1} $. Above, the last inequality is valid because of  $\smash{\frac{k}{2}-\frac{1}{4}<0}$ and $L_{4n}<1.$ 

Finally, we can start with bounding the terms $M_1(z), M_2(z)$, and $M_3$. For $z \in D_1$ and $i \in [n]$,  
we get
\begin{align*}
\left \vert F_i(Z_i(z)) - Z_i(z) + \frac{\rho_i^2}{Z_i(z)} \right \vert &= \left \vert  - \frac{r_{n,i}(z)}{ 1+ r_{n,i}(z)}Z_i(z) + \frac{\rho_i^2}{Z_i(z)}  \right \vert \\  &\leq \frac{10}{9} \left \vert r_{n,i}(z)Z_i(z) - \frac{\rho_i^2(1+r_{n,i}(z))}{Z_i(z)}   \right \vert = \frac{10}{9} \frac{\vert s_i(z) \vert}{\vert Z_i^2(z) \vert},
\end{align*}
where $s_i(z)$ is defined by
\begin{align*}
s_i(z) :=  Z_i^2(z) \left( r_{n,i}(z)Z_i(z) - \frac{\rho_i^2(1+r_{n,i}(z))}{Z_i(z)}  \right).
\end{align*}
The expansions in \eqref{subordination times cauchy expansion} imply
\begin{align*}
s_i(z) = m_3(\mu_i) - \frac{\rho_i^4}{Z_i(z)} - \frac{\rho_i^2m_3(\mu_i)}{Z_i^2(z)} + \left( 1 - \frac{\rho_i^2}{Z_i^2(z)}\right) \int_{\mathbb{R}} \frac{u^4}{Z_i(z) - u} \mu_i(du), \qquad z \in D_1, i \in [n].
\end{align*}
In particular, together with the estimates 
\begin{align*}
	\rho_i^4 \leq m_4(\mu_i), \qquad \rho_i^2 m_3(\mu_i) \leq m_4(\mu_i)^{\frac{5}{4}}, \qquad \rho_i^2 m_4(\mu_i) \leq m_4(\mu_i)^{\frac{3}{2}},
\end{align*}
it follows
\begin{align*}
	\vert s_i(z) \vert \leq \beta_3(\mu_i) + \frac{m_4(\mu_i)}{\vert Z_i(z) \vert} + \frac{m_4(\mu_i)^{\frac{5}{4}}}{\vert Z_i^2(z) \vert} + \frac{m_4(\mu_i)}{\Im z}  + \frac{ m_4(\mu_i)^{\frac{3}{2}}}{\Im z \vert Z_i^2(z) \vert} 
\end{align*}
for all $ z \in D_1$ and $
i \in [n].$ Note that 
\begin{align*}
	\vert M_1(z) \vert = \vert  Z_1(z)  \vert  \left \vert\sum_{i=2}^n F_i(Z_i(z)) - Z_i(z) + \frac{\rho_i^2}{Z_i(z)} \right \vert \leq \frac{10}{9} \frac{1}{\vert Z_1(z) \vert} \sum_{i=2}^n \left \vert \frac{Z_1^2(z)}{Z_i^2(z)}\right \vert \vert s_i(z) \vert
\end{align*}
holds for any $z \in D_1.$ Using \cref{subordination equations} and  \eqref{bound r 0.1}, we obtain 
\begin{align*}
\left \vert \frac{Z_1(z)}{Z_i(z)} - 1 \right \vert = \left \vert \frac{Z_1(z)G_1(Z_1(z))}{Z_i(z)G_i(Z_i(z))} - 1 \right \vert = \left \vert \frac{1+r_{n,1}(z)}{1+r_{n,i}(z)} - 1 \right \vert \leq \frac{10}{9} \left( \vert r_{n,1}(z) \vert + \vert r_{n,i}(z) \vert \right) < \frac{2}{9}
\end{align*}
for $z$  and $i$ as before, which yields
\begin{align} \label{estimate Z_1/Z_i -1}
	\left \vert \frac{Z_1(z)}{Z_i(z)}  \right \vert < \frac{11}{9}, \qquad z \in D_1, i \in [n].
\end{align}
By the previously established bound on $\vert s_i(z) \vert$, \eqref{inequality lyapunov}, and \eqref{lower bound subordination}, we conclude
\begin{align*}
\vert M_1(z) \vert &\leq \frac{10}{9} \frac{121}{81} \frac{1}{\vert Z_1(z) \vert} \sum_{i=2}^{n} \vert s_i(z) \vert \\ & \leq \frac{10}{9} \frac{121}{81} \frac{1}{C_2L_{4n}^{\frac{1}{4} - \frac{k}{2}}} \left( L_{4n}^{\frac{1}{2}} + \frac{L_{4n}}{C_2L_{4n}^{\frac{1}{4} - \frac{k}{2}}} + \frac{L_{4n}^{\frac{5}{4}}}{C_2^2 L_{4n}^{\frac{1}{2} -k}} + \frac{L_{4n}}{\Im z} + \frac{L_{4n}^{\frac{3}{2}}}{\Im z C_2^2 L_{4n}^{\frac{1}{2} -k}}\right) \\ & \leq \frac{10}{9} \frac{121}{81} \frac{1}{C_2} \left( L_{4n}^{\frac{1}{4} + \frac{k}{2}} + \frac{1}{C_2}L_{4n}^{\frac{1}{2} + k} + \frac{1}{C_2^2} L_{4n}^{\frac{1}{2} + \frac{3k}{2}} + \frac{1}{C_1}L_{4n}^{\frac{1}{2}} + \frac{1}{C_1 C_2^2}L_{4n}^{\frac{1}{2} + k} \right) \\ & \leq C_3L_{4n}^{\frac{1}{4} + \frac{k}{2}}
\end{align*}
for all $z \in D_1$ with $C_3 = C_3(k)>0$ being some suitably chosen constant.

Let us continue by bounding the term $M_2(z)$. For this purpose, we combine \eqref{subordination times cauchy expansion} and  \eqref{estimate Z_1/Z_i -1} with the inequality $\rho_1^2 \leq \rho_i^2$, $i \in [n]$, leading to 
\begin{align*}
\left \vert \frac{Z_1(z)}{Z_i(z)} - 1 \right \vert \leq \frac{10}{9} \left( \vert r_{n,1}(z) \vert + \vert r_{n,i}(z) \vert \right) &\leq
\frac{10}{9} \frac{1}{\Im z} \left( \frac{\rho_1^2}{\vert Z_1(z) \vert} + \frac{\rho_i^2}{\vert Z_i(z) \vert}\right) \\ & \leq \frac{10}{9} \frac{1}{\Im z} \frac{\rho_i^2}{\vert Z_1(z) \vert} \left(  1 + \frac{\vert Z_1(z) \vert}{\vert Z_i(z) \vert}\right) < 3 \frac{\rho_i^2}{\Im z \vert Z_1(z) \vert}
\end{align*}
for any $z \in D_1$ and $i \in [n].$
Finally, by making use of \eqref{lower bound subordination}, it follows
\begin{align*}
\vert M_2(z) \vert = \vert Z_1(z) \vert \left \vert \sum_{i=2}^{n} \frac{\rho_i^2}{Z_1(z)} - \frac{\rho_i^2}{Z_i(z)} \right \vert \leq \sum_{i=2}^{n} \rho_i^2 \left \vert 1 - \frac{Z_1(z)}{Z_i(z)} \right \vert < \frac{3}{\Im z \vert Z_1(z) \vert} \sum_{i=2}^{n} \rho_i^4 \leq \frac{3}{C_1 C_2} L_{4n}^{\frac{1}{2}}
\end{align*}
for all $z \in D_1.$

Due to $\rho_1^2 = \rho_1^2 \left( \sum_{i=1}^{n} \rho_i^2 \right)^2 \leq  \left(\sum_{i=1}^{n} \vert \rho_i \vert^3 \right)^2 \leq L_{3n}^2 \leq L_{4n} $, the term $M_3$ admits the estimate 
\begin{align} \label{bound sigma_1}
\vert M_3 \vert = \rho_1^2 \leq L_{4n}.
\end{align}
Consequently, for any $z \in D_1$, we get
 \begin{align} \label{bound q}
	 \vert q(z) \vert = \vert M_1(z) + M_2(z) + M_3 \vert \leq C_3 L_{4n}^{\frac{1}{4} + \frac{k}{2}} + \frac{3}{C_1 C_2} L_{4n}^{\frac{1}{2}} + L_{4n} < C_4  L_{4n}^{\frac{1}{4} + \frac{k}{2}} 
\end{align}
for some constant $C_4 = C_4(k)>0.$

\subsubsection*{Step 2: Analyzing the roots of $Q$}
In this step, we study the roots of $Q(z, \cdot)$ for $z$ in some subset $D_1'$ of $D_1$. Since $Q(z, Z_1(z)) = 0$ holds for all $z \in \mathbb{C}^+$, this will provide information about the behavior of the subordination function $Z_1$ in $D_1'.$

We define
\begin{align*}
D_1':= \left \{	 z \in \mathbb{C}^+ :  \Im z \geq C_5 L_{4n}^{\frac{1}{4} + \frac{k}{2}} \right\} \subset D_1, \qquad C_5 = C_5(k):= \max\{1, C_1, 7C_4\}.
\end{align*}
For the rest of this proof, suppose that $L_{4n}<(3C_5)^{-4}$ is satisfied. We will remove this assumption at the end of the sixth step. 

Now, fix $z \in D_1'$ and let $\omega_1(z), \omega_2(z)$ denote the roots of $Q(z, \cdot).$ Clearly, we have
\begin{align} \label{roots formula}
\omega_i(z) = \frac{1}{2} \left( z + (-1)^i \sqrt{z^2-4+4q(z)}\right), \qquad  i=1,2.
\end{align}
We claim that $\omega_1(z) \neq \omega_2(z)$ holds. Assuming the contrary, we obtain $z^2-4+4q(z) = 0$. Together with $Q(z, Z_1(z)) = 0$, it follows
\begin{align*}
Z_1^2(z) - zZ_1(z)  = - 1 + q(z) = - \frac{z^2 }{4}, 
\end{align*}
which implies $Z_1(z) = \frac{z}{2}$. Due to $\Im Z_1(z) \geq \Im z$ (see \cref{subordination equations}), we arrive at a contradiction. Hence, we get $\omega_1(z) \neq \omega_2(z)$ as claimed.

Let us continue by proving the identity $Z_1(z) = \omega_2(z)$. Suppose that we have $Z_1(z) = \omega_1(z)$. Then,
 \cref{subordination equations} and \eqref{roots formula} yield
\begin{align*}
\Im z \leq \Im \omega_1(z) = \frac{1}{2} \Im z - \frac{1}{2} \Im \sqrt{z^2 -4 + 4q(z)}. 
\end{align*}
If $z^2 -4 + 4q(z) \in \mathbb{C} \setminus [0, \infty)$ is satisfied, we get $\Im \sqrt{z^2 -4 + 4q(z)} \geq 0$ (see \eqref{square root real and im formula}), which leads to the contradiction $\Im z \leq \Im \omega_1(z) \leq \frac{\Im z}{2}$. Thus, it remains to verify $z^2 -4 + 4q(z) \in \mathbb{C} \setminus [0, \infty)$. For this, assume $\Im (z^2 -4 + 4q(z)) = 0$. We obtain $2 \Re z \Im z + 4\Im q(z) = 0$, from which we derive  $	\Re z = - 2\Im q(z)(\Im z)^{-1}.$
Using \eqref{bound q}, the definition of $C_5$, and $L_{4n}< (3C_5)^{-4}$, it follows
\begin{align*}
\Re (z^2 - 4 + 4q(z) )= (\Re z)^2 -  (\Im z)^2 - 4+ 4\Re q(z) &< (\Re z)^2  - 4+ 4\Re q(z) 
= \frac{4 (\Im q(z))^2}{(\Im z)^2} - 4 + 4 \Re q(z) \\ & \leq  \frac{4 \vert q(z) \vert^2}{(\Im z)^2} - 4 +4 \vert q(z) \vert \leq \frac{4C_4^2}{C_5^2} -4 + 4C_4L_{4n}^{\frac{1}{4 }+ \frac{k}{2}} \\ & \leq -3 + 4C_4L_{4n}^{\frac{1}{4 }+ \frac{k}{2}} < -2 < 0.
\end{align*}
Consequently, we arrive at $z^2 - 4 + 4 q(z)\in \mathbb{C} \setminus [0, \infty)$ and $Z_1(z)= \omega_2(z)$ for any $z \in D_1'.$ By similar methods, one can prove $z^2-4 \in \mathbb{C} \setminus [0, \infty)$ for all $z \in D_1'$.
\subsubsection*{Step 3: Application of Bai's inequality}
Later, we will use Bai's inequality from \cref{Bai} with the following parameters: 
\begin{align} \label{choices Bai}
a = 2, \qquad \gamma>0.7, \qquad v = C_5L_{4n}^{\frac{1}{4} + \frac{k}{2}} \in (0,1), \qquad \varepsilon= 6v \in (2va, 2).
\end{align}	
Note that the premise of \cref{Bai} is satisfied due to $m_2(\mu_{S_n}) < \infty.$ Finally, it remains to bound the integrals
\begin{align} \label{integral wrt Re}
 \int_{-\infty}^{\infty} \left\vert G_{S_n}(u+i) - G_\omega(u+i) \right\vert du
 \end{align}
 and 
 \begin{align} \label{integral wrt Im}
 \sup_{x \in I_\varepsilon} \int_{v}^1 \vert G_{S_n}(x+iy) - G_\omega(x+iy) \vert dy
\end{align}
for $I_\varepsilon := [-2+\frac{\varepsilon}{2}, 2-\frac{\varepsilon}{2}] $. 

\subsubsection*{Step 4: Bounding the integral in \eqref{integral wrt Im}} 
In order to bound the integral in \eqref{integral wrt Im}, we have to derive an appropriate estimate for the integrand $\vert G_{S_n}(z) - G_\omega(z) \vert$ for certain $z \in \mathbb{C}^+$. To this end, we write
\begin{align} \label{CT triangle inequality}
	\vert G_{S_n}(z) - G_\omega(z) \vert \leq \left \vert G_{S_n}(z) - \frac{1}{Z_1(z)} \right \vert + \left \vert \frac{1}{Z_1(z)} - G_{\omega}(z) \right \vert, \qquad z \in \mathbb{C}^+,
\end{align}
and study the contributions of the two summands above separately. 

We begin with the second summand. Define $S(z) := F_\omega(z) = \frac{1}{G_\omega(z)}$ and recall that
\begin{align*}
S(z) = \frac{1}{2}\left(z + \sqrt{z^2-4}\right), \qquad z \in \mathbb{C}^+,
\end{align*}
holds. Together with $Z_1= \omega_2$ in $D_1'$  and \eqref{roots formula}, it follows
	\begin{align} \label{difference 1/Z_1 - 1/S}
\frac{1}{Z_1(z)} - G_\omega(z) = 	\frac{1}{Z_1(z)}  - \frac{1}{S(z)}=  \frac{1}{Z_1(z)S(z)} \frac{-2q(z)}{\sqrt{z^2-4} + \sqrt{z^2-4+4q(z)}}
\end{align}
for any $z \in D_1'.$  Let us continue by proving that the real parts of the square roots above have the same sign for all $z \in D_1'$ with $\vert \Re z \vert \geq 1$. The first identity in \eqref{square root real and im formula} yields
\begin{align*}
\sgn \left(\Re \sqrt{z^2-4}\right)  = \sgn(2 \Re  z \Im z) = \sgn(\Re z), \qquad \sgn \left( \Re  \sqrt{z^2-4+4q(z)}\right) = \sgn(2 \Re z \Im z + 4\Im q(z))
\end{align*}
for all $z \in D_1'.$
Now, for $z \in D_1'$ with $\Re z \geq 1$, the definition of $C_5$ combined with \eqref{bound q} implies
 \begin{align*}
2 \Re z \Im z + 4\Im q(z) \geq 2\Im z -  4\vert q(z) \vert > 2C_5 L_{4n}^{\frac{1}{4}+\frac{k}{2}} - 4C_4L_{4n}^{\frac{1}{4} + \frac{k}{2}} > 0,
\end{align*}
whereas we obtain 
\begin{align*}
	2 \Re z \Im z + 4\Im q(z) \leq -2\Im z +  4\vert q(z) \vert < -2C_5 L_{4n}^{\frac{1}{4}+\frac{k}{2}} + 4C_4L_{4n}^{\frac{1}{4} + \frac{k}{2}} < 0
\end{align*}
for $z \in D_1'$ with $\Re z \leq -1$. Thus, as claimed, we have 
\begin{align*}
	\sgn \left(\Re \sqrt{z^2-4}\right)  = \sgn \left( \Re \sqrt{z^2-4+4q(z)}\right)
\end{align*}
for all $z \in D_1'$ with $\vert \Re z \vert \geq 1.$ Since both square roots have positive imaginary part (see \eqref{square root real and im formula}), we derive
\begin{align*}
\left \vert \sqrt{z^2-4} + \sqrt{z^2-4+4q(z)} \right \vert \geq \left \vert \sqrt{z^2-4}  \right \vert, \qquad z \in D_1', \vert \Re z \vert \geq 1.
\end{align*}
Together with \eqref{bound CT omega 1}, \eqref{bound q}, and the inequality
\begin{align*}
\left \vert z^2 -4 \right \vert \geq \max\big\{ \Im z, \big((\Re z)^2-5\big)_+\big\}
\end{align*}
holding for all $z \in \mathbb{C}^+$, it follows
	\begin{align}  \label{estimate 1}
\begin{split}
\left \vert 	\frac{1}{Z_1(z)}  - \frac{1}{S(z)} \right \vert & \leq  \frac{1} {\vert Z_1(z) \vert \vert S(z) \vert} \frac{2 \vert q(z) \vert}{\vert \sqrt{z^2-4} \vert} \\ &\leq \frac{1}{\vert Z_1(z) \vert} \frac{2C_4 L_{4n}^{\frac{1}{4} + \frac{k}{2}}}{ \sqrt{\max\{ \Im z, ((\Re z)^2-5)_+\}}} \leq \frac{1}{\vert Z_1(z) \vert} \frac{2C_4 L_{4n}^{\frac{1}{4} + \frac{k}{2}}}{ \sqrt{\Im z}}
\end{split}	
\end{align}
for all $z \in D_1'$ with $\vert \Re z \vert \geq 1.$

Now, let us consider $z \in D_1'$ with $\vert \Re z \vert <1.$ For such $z$, we have $(\Re z)^2-(\Im z)^2-2  \leq 0$. A simple calculation shows that the last inequality implies
\begin{align*}
	\left(\Im \sqrt{z^2-4} \right)^2 = \frac{1}{2}\left( \sqrt{((\Re z)^2 - (\Im z)^2 -4)^2 + 4(\Re z)^2(\Im z)^2} -\left((\Re z)^2 - (\Im z)^2 -4 \right)\right) \geq 1,
\end{align*}
which in turn leads to 
\begin{align*}
		\left 	\vert \sqrt{z^2-4} + \sqrt{z^2-4+4q(z)}  \right \vert \geq  \Im \left( \sqrt{z^2-4} + \sqrt{z^2-4+4q(z)}\right) \geq  \Im \sqrt{z^2-4} \geq 1.
\end{align*}
By  \eqref{bound CT omega 1}, \eqref{bound q}, and \eqref{difference 1/Z_1 - 1/S}, we deduce
	\begin{align}  \label{estimate 2}
	\left \vert 	\frac{1}{Z_1(z)}  - \frac{1}{S(z)} \right \vert \leq  \frac{2 \vert q(z) \vert}{ \vert Z_1(z) \vert \vert S(z) \vert} \leq \frac{2C_4 L_{4n}^{\frac{1}{4} + \frac{k}{2}}}{\vert Z_1(z) \vert} 
\end{align}
for any $z \in D_1'$ with $\vert \Re z \vert < 1$.

We proceed by proving 
\begin{align} \label{lower bound Z_1 constant}
	\vert Z_1(z) \vert > \frac{1}{10}
\end{align}
for all $z \in D_2$ with 
\begin{align*}
	D_2 := \left \{	 z \in \mathbb{C}^+ : \vert \Re z \vert \leq 2, \, 1 \geq \Im z \geq C_5 L_{4n}^{\frac{1}{4} + \frac{k}{2}} \right\} \subset D_1'.
\end{align*}
Assume the contrary, i.e., we have $\vert Z_1(z) \vert \leq \frac{1}{10}$ for some $z \in D_2$. Then, it follows $\vert Z_1(z) - z \vert \leq \frac{1}{10} + \sqrt{5}<3$. Using the identity $Q(z, Z_1(z)) = 0$ and \eqref{bound q}, we arrive at the following contradiction:
\begin{align*}
	\frac{1}{10} \geq \vert Z_1(z) \vert > \frac{1}{3} \vert Z_1(z) \vert \vert Z_1(z) -z \vert = \frac{1}{3} \vert -\!1 + q(z)  \vert \geq \frac{1}{3} \left( 1- C_4L_{4n}^{\frac{1}{4} + \frac{k}{2}} \right) \geq \frac{1}{3}\left(1-\frac{C_4}{3C_5}\right) > \frac{3}{10}.
\end{align*}
Hence, by \eqref{estimate 1}, \eqref{estimate 2}, and \eqref{lower bound Z_1 constant}, we conclude 
\begin{align} \label{estimate 3}
	\left \vert 	\frac{1}{Z_1(z)}  -G_\omega(z) \right \vert  =	\left \vert 	\frac{1}{Z_1(z)}  - \frac{1}{S(z)} \right \vert \leq \begin{cases}
	20C_4 L_{4n}^{\frac{1}{4} + \frac{k}{2}}(\Im z)^{-\frac{1}{2}} & z \in D_2, \vert \Re z \vert \geq 1, \\
		20C_4 L_{4n}^{\frac{1}{4} + \frac{k}{2}} &  z \in D_2, \vert \Re z \vert < 1.
	\end{cases}
\end{align}

It remains to bound the contribution of the first summand in \eqref{CT triangle inequality}. Making use of \cref{subordination equations} and the definition of $r_{n,1}$, we get 
\begin{align} \label{G_Sn - 1/Z_1}
	 G_{S_n}(z) - \frac{1}{Z_1(z)} =  G_1(Z_1(z)) - \frac{1}{Z_1(z)} =  \frac{r_{n,1}(z)}{Z_1(z)}
\end{align}
for all $z \in \mathbb{C}^+$. In particular, together with  \eqref{subordination times cauchy expansion}, \eqref{bound sigma_1}, \eqref{lower bound Z_1 constant}, and $k<\frac{1}{2}$, it follows
\begin{align} \label{estimate 4}
	\left \vert G_{S_n}(z) - \frac{1}{Z_1(z)} \right \vert  \leq \frac{1}{\vert Z_1(z) \vert}  \frac{\rho_1^2}{\Im z \vert Z_1(z) \vert} \leq \frac{100L_{4n}}{\Im z} \leq \frac{100}{C_5} L_{4n}^{\frac{3}{4}-\frac{k}{2}} \leq \frac{100}{C_5} L_{4n}^{\frac{1}{4}+\frac{k}{2}}
\end{align}
for all $z \in D_2$. 

Recalling the choice of $v$ in \eqref{choices Bai} and combining \eqref{estimate 3} with \eqref{estimate 4}, we derive 
 \begin{align*}
\int_{v}^1 \vert G_{S_n}(x+iy) - G_\omega(x+iy) \vert dy \leq \int_{v}^1 \left( \frac{20C_4}{\sqrt{y}}L_{4n}^{\frac{1}{4} + \frac{k}{2}} + \frac{100}{C_5} L_{4n}^{\frac{1}{4}+\frac{k}{2}} \right) dy \leq \left( 40C_4 +  \frac{100}{C_5} \right)L_{4n}^{\frac{1}{4} + \frac{k}{2}}
\end{align*}
for any $x \in [-2,2]$ with $\vert x \vert \geq 1$. Similarly, for $x \in [-2,2]$ with $\vert x \vert < 1$, we obtain
 \begin{align*}
	\int_{v}^1 \vert G_{S_n}(x+iy) - G_\omega(x+iy) \vert dy \leq \int_{v}^1 \left( 20C_4 L_{4n}^{\frac{1}{4} + \frac{k}{2}}+ \frac{100}{C_5} L_{4n}^{\frac{1}{4}+\frac{k}{2}} \right) dy \leq \left( 20C_4 +  \frac{100}{C_5} \right)L_{4n}^{\frac{1}{4} + \frac{k}{2}}.
\end{align*}
Finally, for $I_\varepsilon = [-2+\frac{\varepsilon}{2}, 2-\frac{\varepsilon}{2}]$, it follows
\begin{align} \label{estimate 6}
\begin{split}
 \sup_{x \in I_\varepsilon} \int_{v}^1 \vert G_{S_n}(x+iy) - G_\omega(x+iy) \vert dy &\leq  \sup_{x \in [-2,2]} \int_{v}^1 \vert G_{S_n}(x+iy) - G_\omega(x+iy) \vert dy \\ & \leq \left( 40C_4 +  \frac{100}{C_5} \right)L_{4n}^{\frac{1}{4} + \frac{k}{2}}.
\end{split}	
\end{align}

\subsubsection*{Step 5: Bounding the integral in \eqref{integral wrt Re}}
We proceed similarly to the last step. Due to $L_{4n}<(3C_5)^{-4}$,  we have $ \{ u + i: u \in \mathbb{R}\} \subset D_1'$. In particular, the inequalities in \eqref{estimate 1} are valid for any complex number $u+i$ with $\vert u \vert \geq 1$, and we deduce
\begin{align} \label{difference 1/Z-1 - Gomega, C_1}
\begin{split}
\left \vert 	\frac{1}{Z_1(u+i)}  -G_\omega(u+i) \right \vert  &=  \left \vert 	\frac{1}{Z_1(u+i)}  - \frac{1}{S(u+i)} \right \vert \\ & \leq  \frac{1}{ \vert Z_1(u+i) \vert \vert S(u+i) \vert} \frac{2 \vert q(u+i) \vert}{\vert \sqrt{(u+i)^2-4} \vert} \\ &\leq \frac{1}{ \vert Z_1(u+i) \vert \vert S(u+i) \vert} \frac{2C_4 L_{4n}^{\frac{1}{4} + \frac{k}{2}}}{  \sqrt{\max\{1, (u^2-5)_+\}}}
\end{split}
\end{align}
for all $u \in \mathbb{R}$ with $\vert u \vert \geq 1.$
A simple calculation (see \cite[equation (3.47)]{Neufeld2024a} for details) shows that 
	\begin{align} \label{lower bound S(u+i)}
	\left \vert \frac{1}{S(u+i)}\right \vert =	\left \vert G_\omega(u+i) \right \vert  \leq \frac{2}{ \sqrt{1 + ((\vert u \vert -4)_+)^2}}
\end{align}
holds for all $u \in \mathbb{R}$. Combining the last two estimates with $\vert Z_1 (u+i)\vert \geq 1$ for $u \in \mathbb{R}$ (see \cref{subordination equations}), we get
\begin{align*}
\left \vert 	\frac{1}{Z_1(u+i)}  - G_\omega(u+i) \right \vert  \leq \frac{4C_4 L_{4n}^{\frac{1}{4} + \frac{k}{2}}}{\sqrt{1+((\vert u \vert -4)_+)^2}  \sqrt{\max\{1, (u^2-5)_+\}}},  \qquad u \in \mathbb{R}, \vert u \vert \geq 1.
\end{align*}
Together with 
\begin{align*}
	\left(1+( \vert u \vert-4)^2\right)\left(u^2 -5\right)\geq \frac{1}{100}\left(1+u^2\right)^2
\end{align*}
for any $u \in \mathbb{R}$ with $\vert u \vert \geq 4$, it follows
\begin{align*} 
	\int_{4}^{\infty} \left\vert \frac{1}{Z_1(u+i)} - G_\omega(u+i) \right\vert du \leq  \int_{4}^{\infty} \frac{40C_4L_{4n}^{\frac{1}{4} + \frac{k}{2}}}{1+u^2} du \leq 10C_4L_{4n}^{\frac{1}{4} + \frac{k}{2}} 
\end{align*}
as well as 
\begin{align*} 
	\int_{-\infty}^{-4} \left\vert \frac{1}{Z_1(u+i)} - G_\omega(u+i) \right\vert du \leq 10C_4L_{4n}^{\frac{1}{4} + \frac{k}{2}}.  
\end{align*}
Moreover, we have 
\begin{align*} 
&	\int_{\sqrt{6}}^{4} \left\vert \frac{1}{Z_1(u+i)} - G_\omega(u+i) \right\vert du \leq  \int_{\sqrt{6}}^{4}  \frac{4C_4L_{4n}^{\frac{1}{4} + \frac{k}{2}}  }{\sqrt{u^2-5}}du \leq 4C_4L_{4n}^{\frac{1}{4} + \frac{k}{2}}, \\
	& \qquad \qquad	\int_{-4}^{-\sqrt{6}} \left\vert \frac{1}{Z_1(u+i)} - G_\omega(u+i) \right\vert du  \leq 4C_4L_{4n}^{\frac{1}{4} + \frac{k}{2}},
\end{align*}
and 
\begin{align*}
&	\int_{1}^{\sqrt{6}} \left\vert \frac{1}{Z_1(u+i)} - G_\omega(u+i) \right\vert du \leq  \int_{1}^{\sqrt{6}}  4C_4L_{4n}^{\frac{1}{4} + \frac{k}{2}}  du \leq 6C_4L_{4n}^{\frac{1}{4} + \frac{k}{2}}, \\
& \qquad \qquad \int_{-\sqrt{6}}^{-1} \left\vert \frac{1}{Z_1(u+i)} - G_\omega(u+i) \right\vert du  \leq 6C_4L_{4n}^{\frac{1}{4} + \frac{k}{2}}.
\end{align*}
It remains to fill the gap in the integration from $-1$ to $1$. To this end, we note that 
	\begin{align*}
	\left \vert \sqrt{(u+i)^2 -4} + \sqrt{(u+i)^2 - 4 +4q(u+i)}  \right \vert \geq \Im  \sqrt{(u+i)^2 -4}  \geq 1
\end{align*}
holds for all $u \in \mathbb{R}$. Together with \eqref{difference 1/Z_1 - 1/S} and $\vert S(u+i) \vert, \vert Z_1(u+i)\vert \geq 1$ for $u \in \mathbb{R}$, we obtain
\begin{align*}
\left \vert 	\frac{1}{Z_1(u+i)}  - G_\omega(u+i) \right \vert  \leq 2C_4 L_{4n}^{\frac{1}{4} + \frac{k}{2}},  \qquad u \in \mathbb{R}.
\end{align*} Integration yields
\begin{align*}
		\int_{-1}^{1} \left\vert \frac{1}{Z_1(u+i)} - G_\omega(u+i) \right\vert du \leq 4C_4L_{4n}^{\frac{1}{4} + \frac{k}{2}}. 
\end{align*}
In particular, we conclude
\begin{align} \label{estimate 5}
	\int_{-\infty}^{\infty} \left\vert \frac{1}{Z_1(u+i)} - G_\omega(u+i) \right\vert du \leq 44C_4L_{4n}^{\frac{1}{4} + \frac{k}{2}}.
\end{align}

Finally, we consider the contribution of the difference between $G_{S_n}$ and $\frac{1}{Z_1}$. For this, we need the inequality  	
\begin{align} \label{bounds for abs Z_1 for integral wrt Re of 1/Z_1 - 1/S}
	\vert Z_1(u+i) \vert \geq 	\frac{1}{5}\left( 1 + (\vert u \vert -4)_+\right)
\end{align}
for all $u \in \mathbb{R}$ with $\vert u \vert \geq 1$, 
which can be proved as follows: By \eqref{difference 1/Z-1 - Gomega, C_1} and \eqref{lower bound S(u+i)}, we get
\begin{align*}
\vert Z_1(u+i) \vert &\geq  \vert S(u+i) \vert - \vert Z_1(u+i) - S(u+i) \vert \\ &\geq \frac{ \sqrt{1 + ((\vert u \vert -4)_+)^2}}{2} - \frac{2C_4 L_{4n}^{\frac{1}{4} + \frac{k}{2}}}{  \sqrt{\max\{1, (u^2-5)_+\}}}
\end{align*} 
for any $u \in \mathbb{R}$ with $\vert u \vert \geq 1.$ Now, for $u \in \mathbb{R}$ with $1 \leq \vert u \vert \leq 4$, the inequality in \eqref{bounds for abs Z_1 for integral wrt Re of 1/Z_1 - 1/S} follows from 
\begin{align*}
 \frac{2C_4 L_{4n}^{\frac{1}{4} + \frac{k}{2}}}{  \sqrt{\max\{1, (u^2-5)_+\}}} \leq  2C_4 L_{4n}^{\frac{1}{4} + \frac{k}{2}} \leq \frac{2C_4}{3C_5} < \frac{1}{10}.
 \end{align*}
 For $u \in \mathbb{R}$ with $\vert u \vert >4$, we use
 \begin{align*}
 \sqrt{1 + (\vert u \vert -4)^2} \geq \frac{1}{\sqrt{2}}(1 + (\vert u \vert - 4)) \qquad \text{and} \qquad \frac{2C_4 L_{4n}^{\frac{1}{4} + \frac{k}{2}}}{\sqrt{u^2-5}} < \frac{1}{10},
 \end{align*}
and arrive at
 \begin{align*}
 \frac{ \sqrt{1 + (\vert u \vert -4)^2}}{2} - \frac{2C_4 L_{4n}^{\frac{1}{4} + \frac{k}{2}}}{  \sqrt{ u^2-5}} \geq \frac{1}{\sqrt{2}} (1+ ( \vert u \vert - 4)) - \frac{1}{10} > \frac{1}{5}(1+ (\vert u \vert -4)).
 \end{align*}
Together with  \eqref{G_Sn - 1/Z_1}, \eqref{subordination times cauchy expansion}, and \eqref{bound sigma_1}, we obtain
 \begin{align*}
 \int_{4}^\infty  \left\vert  G_{S_n}(u+i) - \frac{1}{Z_1(u+i)}  \right\vert du &=  \int_{4}^\infty  \left\vert \frac{r_{n,1}(u+i)}{Z_1(u+i)}  \right\vert du \\ &\leq \int_{4}^\infty \frac{\rho_1^2}{\vert Z_1^2(u+i) \vert} du \leq 25 L_{4n} \int_{4}^\infty \frac{1}{(u-3)^2} du = 25L_{4n}
 \end{align*}
 and \begin{align*}
 \int_{-\infty}^{-4} \left\vert G_{S_n}(u+i) - \frac{1}{Z_1(u+i)}  \right\vert du \leq 25L_{4n}.
 \end{align*}
Due to
 \begin{align*}
  \int_{-4}^{4}  \left\vert G_{S_n}(u+i) - \frac{1}{Z_1(u+i)}  \right\vert du = \int_{-4}^4  \left\vert \frac{r_{n,1}(u+i)}{Z_1(u+i)}  \right\vert du \leq 8\rho_1^2 \leq 8L_{4n},
 \end{align*}
 we conclude
 \begin{align*}
\int_{-\infty}^\infty \left\vert G_{S_n}(u+i) - \frac{1}{Z_1(u+i)} \right\vert du  \leq 58L_{4n}.
 \end{align*}
Combining the last inequality with \eqref{estimate 5}, it follows
 \begin{align} \label{estimate 7}
 \int_{-\infty}^\infty \left\vert G_{S_n}(u+i) - G_{\omega}(u+i) \right\vert du  \leq  44C_4L_{4n}^{\frac{1}{4} + \frac{k}{2}} + 58L_{4n} \leq \left( 44C_4 + 58\right)L_{4n}^{\frac{1}{4} + \frac{k}{2}}.
 \end{align}
 
 \subsubsection*{Step 6: Final conclusion} 
For completeness, we end the proof of \cref{main theorem 2} with a summary. First, assume that $L_{4n}<(3C_5)^{-4}$ holds. Then, using \cref{Bai} with the parameters in \eqref{choices Bai} as well as  \eqref{estimate 6} and \eqref{estimate 7}, we obtain 
\begin{align*}
	\Delta(\mu_{S_n}, \omega)  \leq C_\gamma \! \left( \frac{16C_5}{\pi} L_{4n}^{\frac{1}{4} + \frac{k}{2}} +  \gamma (6C_5)^{\frac{3}{2}}L_{4n}^{\frac{3}{8} +\frac{3k}{4}} + \left( 40C_4 +  \frac{100}{C_5} \right)\!L_{4n}^{\frac{1}{4} + \frac{k}{2}} +\left( 44C_4 + 58\right)\!L_{4n}^{\frac{1}{4} + \frac{k}{2}} \right)  \!\leq C_6 L_{4n}^{\frac{1}{4} + \frac{k}{2}}
\end{align*}
for some suitably chosen constant $C_6 = C_6(k) >0.$ 

Second, suppose that we have $L_{4n} \geq (3C_5)^{-4}$. Due to $k<\frac{1}{2}$, it follows
\begin{align*}
	\Delta(\mu_{S_n}, \omega)  \leq 1 \leq 9C_5^2L_{4n}^{\frac{1}{2}} \leq 9C_5^2  L_{4n}^{\frac{1}{4} + \frac{k}{2}}. 
\end{align*}
In particular, by choosing $C(k):= \max\{C_6, 9C_5^2\}$, \cref{main theorem 2} is proven. By definition, $C(k)$ depends on $k$ only via $D(k)$.

  \end{document}